\documentclass[12pt]{article}%
\usepackage[applemac]{inputenc}
\usepackage{amsmath,amssymb}
\usepackage{amsfonts}
\usepackage[applemac]{inputenc}
\usepackage{amsmath,amssymb}
\usepackage{color}
\usepackage{color, colortbl}
\usepackage{amsmath}
\usepackage{amssymb}
\usepackage{graphicx}
\usepackage{tikz}
\usepackage{geometry}
\usetikzlibrary{arrows}
\usepackage{amsfonts}
\usepackage{amsmath,amssymb}
\usepackage[applemac]{inputenc}
\usepackage{color}
\usepackage{amsmath}
\usepackage{amssymb}
\usepackage{graphicx}
\usepackage{tikz}
\newtheorem{theorem}{Theorem}[section]

\newtheorem{lemma}[theorem]{Lemma}
\newcommand{\E}{E}

\newtheorem{definition}[theorem]{Definition}

\newtheorem{example}[theorem]{Example}

\newtheorem{problem}[theorem]{Problem}

\newtheorem{remark}[theorem]{Remark}

\usetikzlibrary{arrows,shapes,automata,petri}
\newenvironment{myenumerate}{

\begin{enumerate}}{\end{enumerate}}
\newcommand{\dproof}{\noindent {Proof.} \quad}
\newcommand{\fproof}{\hfill $\square$ \bigskip}

\numberwithin{equation}{section}

\definecolor{LightCyan}{rgb}{0.88,1,1}

\def\RR{{\mathbb{ R}}}

\def\1B{\text{1\!\!I}}

\def\l{\langle}
\def\<{\langle}
\def\>{\rangle}
\def\P{\mathbb{P}}
\def\R{\mathbb{R}}
\def\l{\lambda}
\def\S{\mathcal{S}}
\def\E{\mathbb{E}}

\def\L{\mathcal{L}}

\begin{document}

\title{The stochastic fractional transport heat equation}

\author{Olfa Draouil$^{1}$, Rahma Yasmina Moulay Hachemi$^{2}$ \& Bernt \O ksendal$^{3}$
}
\date{2 September 2026}
\vskip 0.5cm

\maketitle

\footnotetext[1]{Department of Mathematics, University of Tunis El Manar, Tunisia\newline
Email: olfa.draouil@fst.utm.tn} 

\footnotetext[2]{%
Department of Mathematics, University of Oslo, Norway.\newline
Email: yasmin.moulayhachemi@yahoo.com}

\footnotetext[3]{%
Department of Mathematics, University of Oslo, Norway. 
Email: oksendal@math.uio.no.}


\begin{abstract}
We give an introduction to the  time-fractional stochastic heat equation driven by 1+d-parameter fractional time-space white noise, in the following two cases:
\begin{myenumerate}
    \item
With additive fractional time-space white noise
  \item 
With multiplicative time-space Brownian noise
  \end{myenumerate}
  The fractional time derivative is interpreted as the Caputo derivative of order $\alpha \in (0,2]$ and we assume that the Hurst coefficient $H=(H_0,H_1,H_2, ...,H_d)$ of the time-space fractional white noise is in $(\tfrac{1}{2},1)^{1+d}$.

  We find an explicit expression for the unique solution in the sense of distribution of the equation in the additive noise case (i).
  
  In the multiplicative case (ii) we derive the correct space--time convolution equation and its Wiener chaos expansion. For $1<\alpha\leq2$ a second initial condition is required. The endpoint $\alpha=2$ is treated separately and yields a stochastic wave equation with finite propagation speed.

A solution $Y(t,x)$ is called \emph{mild} if $E[Y^2(t,x)] < \infty$ for all $t,x$. In the multiplicative case the random-field condition is expressed by the square integrability of the forcing kernel. For the classical Laplacian and space--time white noise, the endpoint $\alpha=2$ admits a random-field solution only in space dimension $d=1$.

This paper is partly a survey paper, explaining the concepts and methods behind the results. It is also partly a research paper, in the sense that some results are new, to the best of our knowledge.

  \end{abstract}

  \section{Introduction}
\indent
The purpose of this paper is to study the time-fractional heat equation driven by multi-parameter fractional Brownian motion. We discuss both the additive noise case and the case with multiplicative noise. Using basic methods of Laplace and Fourier transforms combined with white noise calculus,  we refer to \cite{streit}, we obtain explicit solution formulas in the sense of distribution for any order $\alpha \in (0,2)$ of the time -derivative and for any space dimension $d$. Then we investigate for what values of $\alpha$ and $d$ the solution has finite second moment.

It was Niels Henrik Abel who first introduced the fractional derivative of a function.
In 1823 [1] he used the fractional derivative to solve the tautochrone (isochrone)
problem in mechanics.
Then in 1903 [22]  Gosta Magnus Mittag-Leffler introduced what was later called the the Mittag-Leffler function $E_{\alpha}(z)$. Subsequently it became clear that this function has a connection to
the fractional derivative introduced by Abel.

As a modelling tool the fractional derivative of order $\alpha \in (0,2)$ is useful for many situations, e.g. in the
study of heat transfer and waves. For many applications
of fractional derivatives we refer to the book by S. Holm [15].

 - The \emph{classical diffusion case} corresponds to $\alpha=1$. In that case our equation models the normal diffusion of
heat in a random or noisy medium, the noise being represented by the fractional
time-space white noise $W_H (t, x)$.

- The \emph{diffusion--wave case} corresponds to $1<\alpha<2$. It interpolates between diffusion and wave propagation and requires both an initial value and an initial time derivative. At $\alpha=2$ the equation becomes a wave equation. Fractional diffusion--wave models may describe transport with memory and propagating fronts.

- The \emph{subdiffusion case}  occures when $\alpha <1$. Then the equation models a situation in which travel times of the
particles are longer than in the standard case. Such situations may occur for example in some transport systems. For more details see \cite{cherstry}.

Specifically, in this paper we study the  time-fractional stochastic heat equation in the following two cases:
\begin{myenumerate}
\item
Additive noise:
\begin{align}\label{heat0}
 	 &\frac{\partial^{\alpha}}{\partial t^{\alpha}}Y(t,x)=-\lambda (-\Delta)^{ \beta} Y(t,x)+b \bullet \nabla Y(t,x)+\sigma W_{H}(t,x)\; (t,x)\in (0,\infty)\times \mathbb{R}^{d}.
	 \end{align}
	 where $W_H$ denotes $1+d$ -parameter  time-space fractional white noise
	 \item
Multiplicative noise:
 \begin{align}\label{heat1}
 	 &\frac{\partial^{\alpha}}{\partial t^{\alpha}}Y(t,x)=-\lambda (-\Delta)^{ \beta} Y(t,x)+b \bullet \nabla Y(t,x)+\sigma Y(t,x)\diamond W(t,x)\; (t,x)\in (0,\infty)\times \mathbb{R}^{d}.
	 \end{align}
where $W$ denotes classical $1+d$ -parameter time-space white noise.
\end{myenumerate}
	  	where $d\in\mathbb{N}=\{1,2,...\}$ and $\frac{\partial^{\alpha}}{\partial t^{\alpha}}$ is the Caputo derivative of order $\alpha \in (0,2]$, and $\l>0$ and $\sigma\in \mathbb{R}$ are given constants. 
Here $\diamond$ denotes the Wick product,
\begin{equation}
  		\Delta Y =\sum_{j=1}^{d}\frac{\partial^{2}Y}{\partial x_{j}^{2}}(t,x)
  	\end{equation} 
denotes the Laplacian operator acting on $x$,  $W_{H}(t,x)$ and $W(t,x)$ denote multiparameter fractional  and classical white noise, respectively, defined by
\begin{equation}
  	W_{H}(t,x)=\frac{\partial}{\partial t}\frac{\partial^{d}B_{H}(t,x)}{\partial x_{1}...\partial x_{d}},
  \end{equation}
\begin{equation}
  	W(t,x)=\frac{\partial}{\partial t}\frac{\partial^{d}B(t,x)}{\partial x_{1}...\partial x_{d}},
  \end{equation}
\noindent where $B(t,x)=B[t,x,\omega); t\geq 0, x \in \R^d, \omega \in \Omega$ is classical $1+d$ -parameter Brownian sheet and $B_H(t,x,\omega); t\geq 0, x \in \R^d, \omega \in \Omega$  is time-space fractional Brownian
 sheet with Hurst parameter $H=(H_0, H_1, ... , H_d) \in (\tfrac{1}{2},1)^{d+1}$.

The boundary conditions are
\begin{align}
    Y(0,x)&=\delta_0(x)\text{ (the point mass at  } 0), \label{1.4}\\
    \partial_tY(0,x)&=0,\qquad 1<\alpha\leq2,\nonumber\\
    \lim_{x \rightarrow +/- \text{ }\infty}Y(t,x)&=0.\label{1.5}
\end{align}
 We consider the equation \eqref{heat0} in the sense of distribution, and we find an explicit expression for the $\mathcal{S}'$-valued solution $Y(t,x)$, where $\mathcal{S}'$ is the space of tempered distributions. See Theorem \ref{add}.\\
 
The equation \eqref{heat1} is considered both in the Hida space $(\mathcal{S})^{*}$ and, when the kernel is square integrable, as a Walsh random-field equation. We derive its mild space--time convolution form in Theorem \ref{mult1} and its Wiener chaos expansion in Theorem \ref{vol3}.

The solution is called \emph{mild} if $Y(t,x) \in L^2(\P)$ for all $t,x$. For the multiplicative equation driven by classical space--time white noise, the exact criterion is the square integrability condition \eqref{multdalang}; see Theorem \ref{mild2}. In particular, when $\alpha=2$ and $\beta=1$, a random-field solution exists only for $d=1$. The distribution-valued formulation remains available beyond the random-field regime.\\
 
 This paper extends the results  in our previous paper \cite{MO} in the following sense:\\
 (i) The additive case is extended from classical Brownian sheet to fractional Brownian motion\\
 (ii) The multiplicatice case is not covered in \cite{MO}.\\
 

There are many papers dealing with various forms of stochastic fractional differential equations. Some papers which are related to ours are:  \\
 -  In the paper by Kochubel et al \cite{KKd} the fractional heat equation corresponding to random time change in Brownian motion is studied.\\
 - The papers by Bock et al \cite{BGO}, \cite{Bd} 
 are considering stochastic equations driven by grey Brownian motion. \\
  - The paper Roeckner et al \cite{LRd} proves existence and uniqueness of general time-fractional linear evolution equations in the Gelfand triple setting.\\
  -The paper Yalin et al \cite{Y} study the time-regularity of the paths of solutions to stochastic
  partial differential equations  driven by additive
  infinite-dimensional fractional Brownian noise.\\
  -The paper by Binh et al  \cite{Bi} study the spatially-temporally Holder continuity of mild random field solutions of space time fractional stochastic heat equation driven by colored noise.\\
  - The paper  which is closest to our paper is Chen et al \cite{CHN}, where 
a comprehensive discussion is given of a general fractional stochastic heat equations with multiplicative noise, and with fractional derivatives in both time and space, is given. In that paper the authors prove existence and uniqueness results as well as regularity results of the solution, and they give sufficient conditions on the coefficients and the space dimension $d$, for the solution to be a random field. Their method is based on the use of fundamental solutions of the corresponding deterministic equations.

In comparison, our paper  is dealing with both the additive and multiplicative noise case. We consider a more special class of fractional heat equations than in Chen et al \cite{CHN} and this allows us to use basic methods of Laplace and Fourier transform and to find explicit solution formulas in the sense of distributions.
Then we use this explicit solution formulas to find conditions under which the solutions are random fields in $L^2(\P).$\\

We refer to Holm \cite{H}, Ibe \cite{I}, Kilbas et al \cite{Kilbas} and Samko et al \cite{Samko} for more information about fractional calculus and their applications.

\section{Preliminaries}
In this section we summarise some of the background notation and results needed in this paper.
\subsection{The space of tempered distributions}
For the convenience of the reader we recall some of the basic properties of the Schwartz space $\mathcal{S}$ of rapidly decreasing smooth functions and its dual, the space $\mathcal{S}'$ of tempered distributions.
Let $n$ be a given natural number. Let $\mathcal{S}=\mathcal{S}(\mathbb{R}^n)$\label{simb-028} be the
space of rapidly decreasing smooth real
functions $f$
on $\mathbb{R}^n$
equipped with the family of seminorms:\label{simb-029} 
\begin{equation*}
\Vert f \Vert_{k,\alpha} := \sup_{y \in \mathbb{R}^n}\big\{ (1+|y|^k) \vert
\partial^\alpha f(y)\vert \big\}< \infty,
\end{equation*}
where $k = 0,1,...$, $\alpha=(\alpha_1,...,\alpha_n)$ is a multi-index with $%
\alpha_j= 0,1,...$ $(j=1,...,n)$ and\label{simb-030} 
\begin{equation*}
\partial^\alpha f := \frac{\partial^{|\alpha|}}{\partial
y_1^{\alpha_1}\cdots \partial y_n^{\alpha_n}}f
\end{equation*}
for $|\alpha|=\alpha_1+ ... +\alpha_n$.

Then
$\mathcal{S}=\mathcal{S}(\mathbb{R}^n)$ is a
Fr\'echet space.

Let $\mathcal{S}^{\prime }=\mathcal{S}^{\prime }(\mathbb{R}^{n})$\label%
{simb-031} be its dual, called the space of \emph{tempered distributions}. 
\index{tempered distributions} Let $\mathcal{B}$ denote the family of all
Borel subsets of $\mathcal{S}^{\prime }(\mathbb{R}^{n})$ equipped with the
weak* topology. If $\Phi \in \mathcal{S}^{\prime }$ and $f \in \mathcal{%
S}$ we let \label{simb-033} 
\begin{equation}
\Phi (f) \text{ or } \langle \Phi ,f \rangle  \label{3.1}
\end{equation}%
denote the action of $\Phi$ on $f$.

\begin{example}
\begin{itemize}
\item
{(Evaluation at a point)}
For $y \in \R$ define the function $\delta_y$ on $\S(\R)$ by $\delta_y(\phi)=\phi(y)$. Then $\delta_y$ is a tempered distribution.\\
\item
{(Derivatives)} Consider the function $D$, defined for $\phi \in \S(\mathbb{R})$ by $D[\phi]=\phi^{\prime}(y)$. Then  $D$ is a tempered distribution. \\
\item
{(Distributional derivative)}\\
 Let $T$ be a tempered distribution, i.e. $T \in \S^{'}(\mathbb{R}) $. We define the distributional derivative $T^{'}$ of $T$ by
 $$ T^{'}[\phi]=-T[\phi^{'}]; \quad \phi \in \S.$$
 Then $T^{'}$ is again a tempered distribution.
 \end{itemize}
 \end{example}
 
 In the following we will apply this to the case when $n=1+d$ and $y=(t,x) \in \R \times \R^d$.
 \subsection{The Mittag-Leffler functions}
 \begin{definition}(The two-parameter Mittag-Leffler function)\newline
 The Mittag-Leffler function of two parameters $\alpha,\beta$ is denoted by $E_{\alpha,\beta}(z)$ and defined by:
\begin{equation}
    E_{\alpha,\beta}(z)=\sum_{k=0}^{\infty}\frac{z^{k}}{\Gamma(\alpha k+\beta)}
\end{equation}
where $z,\; \alpha,\; \beta\in \mathbb{C}$ (the set of complex numbers),\; $Re(\alpha)>0\; and\; Re(\beta)>0,$ and $\Gamma$ is the Gamma function.\\

 For $\beta=1$ we obtain
the Mittag-Leffler function of one parameter $\alpha$ denoted by $E_{\alpha}(z)$. This function is defined as
\begin{equation}
    E_{\alpha}(z)=\sum_{k=0}^{\infty}\frac{z^{k}}{\Gamma(\alpha k+1)}
\end{equation}
where $z,\; \alpha\in \mathbb{C},\; Re(\alpha)>0.$
\end{definition}

\begin{remark} Note that $E_{\alpha}(z)= E_{\alpha,1}(z)$ and that 
\begin{align}
E_{1}(z)=\sum_{k=0}^{\infty}\frac{z^{k}}{\Gamma( k+1)}=\sum_{k=0}^{\infty}\frac{z^{k}}{k!} =e^{z}.
\end{align}
\end{remark}

\subsection{The (Abel-)Caputo fractional derivative}
  In this section we present the definitions and some properties of the Caputo derivatives.
  \begin{definition}
  The (Abel-)Caputo fractional derivative of order $\alpha > 0$ of a function $f$ such that $f(x)=0$ when $x<0$ is denoted by  $D^{\alpha} f (x)$ or $\frac{d^{\alpha}}{dx^{\alpha}} f(x)$ and defined  by 
  \begin{align}\label{caputo1}
  D^{\alpha}f(x):& =
  \begin{cases}
  \frac{1}{\Gamma(n-\alpha)}\int_0^x \frac{f^{(n)}(u)du}{(x-u)^{\alpha +1 -n}}; \quad n-1 < \alpha < n\\
  \frac{d^n}{dx^n}f(x); \quad \alpha =n.
  \end{cases}
  \end{align}
  Here $n$ is an smallest integer greater than or equal to $\alpha$.\\
  
  \noindent If $f$ is not smooth these derivatives are interpreted in the sense of distributions.
  \end{definition}
  
  \begin{example}
  If $f(x)=x$ and $\alpha \in (0,1)$ then
 
  \begin{align}
 D^{\alpha}f(x)=\frac{ x^{1-\alpha} }{(1-\alpha)\Gamma(1-\alpha)}.
 \end{align} 
 In particular, choosing $\alpha=\tfrac{1}{2}$ we get
 \begin{align}
 D^{\tfrac{1}{2}}f(x)=\frac{2 \sqrt{x}}{\sqrt{\pi}}.
 \end{align} 
  \end{example}

\subsection{Laplace transform of Caputo derivatives} 
\begin{itemize}
\item
The 1-dimensional case.\\
Recall that the (1-dimensional) Laplace transform $L$ is defined by 
 
\begin{equation}
	Lf(s)=\int_{0}^{\infty}e^{-st}f(t)dt=\widetilde{f}(s)
	\end{equation}
	for all $f$ such that the integral converges.

 Some of the properties of the Laplace transform that we will need are:

  \begin{align}
    &L[ \frac{\partial ^{\alpha}}{\partial t^{\alpha}}f(t)](s)=s^{\alpha}(L f)(s)-s^{\alpha-1}f(0) \label{L1}\\       
     &L[E_{\alpha}(bx^{\alpha})](s)  = \frac{s^{\alpha -1}}{s^{\alpha}-b}\label{L2}\\
     &L[x^{\alpha-1}E_{\alpha,\alpha}(-b x^{\alpha})](s)=\frac{1}{s^{\alpha}+b}\label{L3}
     \end{align}
     Recall that the convolution $f\ast g$ of two functions $f,g: [0,\infty) \mapsto \mathbb{R}$  is defined by
\begin{align}
(f \ast g)(t)=\int_0^t f(t-r)g(r) dr; \quad t \geq 0.
\end{align}

The convolution rule for Laplace transform states that $$L\left( \int_{0}^{t}f(t-r)g(r)dr\right) (s)=Lf(s)Lg(s),$$ 
or 
\begin{equation}\label{12aa}
	\int_{0}^{t}f(t-r)g(r)dr=L^{-1}\left( Lf(s)Lg(s)\right) (t).
\end{equation}

\item
{The multi-dimensional case}\\
If $g: \mathbb{R}_{+}^d \mapsto \mathbb{R}$ is bounded, we define the \emph{multi-dimensional }Laplace transform $\L$ of $g$ by

\begin{equation}
	\L\{g\}(y)=\int_{\mathbb{R}^d_{+}}e^{-xy}g(x)dx=:\widehat{g}(y);\quad y\in \mathbb{R}_{+}^d; 
	\end{equation}
 where $xy=x_1y_2 + x_2y_2 + ... +x_dy_d$ when $x=(x_1,x_2, ... ,x_d), y=(y_1,y_2, ... ,y_d).$\\
 In this case the convolution $f\ast g$ of two functions $f,g: [0,\infty)^d \mapsto \mathbb{R}$  is defined by
\begin{align}
(f \ast g)(x)=\int_0^x f(x-z)g(z) dz; \quad x \in \mathbb{R}_{+}^d,
\end{align}
where we use the notation
$\int_0^x = \int_0^{x_1} \int_0^{x_2} ... \int_0^{x_d}$ if $x=(x_1,x_2, ... x_d).$\\

The convolution rule for the multi-dimensional Laplace transform states that $$\mathcal{L}\left( \int_{0}^{x}f(x-z)g(z)dz\right) (y)=\L f(y) \L g(y),$$ 
or 
\begin{equation}\label{12A}
	\int_{0}^{x}f(x-z)g(z)dz=\L ^{-1}\left( \L f(y) \L g(y)\right) (x).
\end{equation}

\end{itemize}

 \section{Time-space white noise}
In this section we briefly review the basic notation and results of time-space white noise. 

Let $n$ be a fixed natural number. Later we will set $n= 1 + d$.
Define $\Omega={\cal S}'(\RR^n)$, equipped
with the weak-star topology. This space will be the base of our basic
probability space, which we explain in the following:
\vskip 0.3cm
As events we will use the family $\mathcal{F}={\cal
B}({\cal S}'(\RR^n))$ of Borel subsets of ${\cal S}'(\RR^d)$, and our
probability measure $\P$ is defined by the following result:

\begin{theorem}{\bf (The Bochner--Minlos theorem)}\\
There exists a unique probability measure $\P$ on ${\cal B}({\cal S}'(\RR^n))$
with the following property:
$$\E[e^{i\langle\cdot,\phi\rangle}]:=\int\limits_{\cal S'}e^{i\langle\omega,
\phi\rangle}d\mu(\omega)=e^{-\tfrac{1}{2} \Vert\phi\Vert^2};\quad i=\sqrt{-1}$$
for all $\phi\in{\cal S}(\RR^n)$, where
$\Vert\phi\Vert^2=\Vert\phi\Vert^2_{L^2(\RR^n)},\quad\langle\omega,\phi\rangle=
\omega(\phi)$ is the action of $\omega\in{\cal S}'(\RR^n)$ on
$\phi\in{\cal S}(\RR^n)$ and $\E=\E_{\P}$ denotes the expectation
with respect to  $\P$.
\end{theorem} 
We will call the triplet $({\cal S}'(\RR^n),{\cal
B}({\cal S}'(\RR^n)),\P)$ the {\it  white noise probability
space\/}, and $\P$ is called the {\it white noise probability measure}.

The measure $\P$ is also often called the (normalised) {\it Gaussian
measure\/} on ${\cal S}'(\RR^n)$. It is not difficult to prove that if $\phi\in L^2(\RR^n)$ and
we choose $\phi_k\in{\cal S}(\RR^n)$ such that $\phi_k\to\phi$ in $L^2(\RR^n)$,
then
$$\langle\omega,\phi\rangle:=\lim\limits_{k\to\infty}\langle\omega,\phi_k\rangle
\quad\text{exists in}\quad L^2(\P)$$
and is independent of the choice of $\{\phi_k\}$. In
particular, if we define
$$\widetilde{B}(x):=\widetilde{B}(x_1,\cdots,x_n,\omega)=\langle\omega,\chi_
{[0,x_1]\times\cdots\times[0,x_n]}\rangle; \quad  
x=(x_1,\cdots,x_n)\in\RR^n,$$
where $[0,x_i]$ is interpreted as $[x_i,0]$ if $x_i<0$,
then $\widetilde{B}(x,\omega)$ has an $x$-continuous version $B(x,\omega)$, which
becomes an \emph{$n$-parameter Brownian motion}, in the following sense:

By an \emph{$n$-parameter Brownian motion} we mean a family
$\{B(x,\cdot)\}_{x\in\RR^n}$ of random variables on a probability space
$(\Omega,{\cal F},\P)$ such that
\begin{itemize} 
\item
$B(0,\cdot)=0\quad\text{almost surely with respect to } \P,$
\item
$\{B(x,\omega)\}$ is a continuous and Gaussian stochastic process 
\item
For all $x=(x_1,\cdots,x_n)$,
$y=(y_1,\cdots,y_n)\in\RR_+^n$,
$B(x,\cdot),\,B(y,\cdot)$ have the covariance
$\prod_{i=1}^n x_i\wedge y_i$. 
\end{itemize} 

It can be proved that the process $\widetilde{B}(x,\omega)$ defined above has a modification $B(x,\omega)$ which satisfies all these properties.
This process $B(x,\omega)$ then becomes an \emph{$n$-parameter Brownian motion}. 
\vskip0.3cm
We remark that for $n=1$ we get the classical (1-parameter) Brownian motion
$B(t)$ if we restrict ourselves to $t\geq 0$. For $n \geq 2$ we get what is often
called the \emph{Brownian sheet}.
\vskip0.3cm
With this definition of Brownian motion it is natural to define the
$n$-parameter Wiener--It\^{o} integral of $\phi\in L^2(\RR^n)$ by
$$\int\limits_{\RR^n}\phi(x)dB(x,\omega):=\langle\omega,\phi\rangle;\quad
\omega\in{\cal S}'(\RR^n).$$
We see that by using the Bochner--Minlos theorem we have obtained an easy construction of
$n$-parameter Brownian motion that works for any parameter dimension $n$. Moreover, we get a representation of the space $\Omega$ as the Fr\' echet space $\mathcal{S}'(\R^n)$. This is an advantage in many situations, for example in the construction of the Hida-Malliavin derivative, which can be regarded as a stochastic gradient on $\Omega$.
See e.g. \cite{Dop} and the references therein.

In the following we put $n=1+d$ and let 
 $$B(t,x)=B(t,x,\omega); t \geq 0, x \in \R^d, \omega \in \Omega$$ 
 denote the (1-dimensional) time-space Brownian motion (also called the Brownian sheet) with probability law $\P$. Since this process is $(t,x)$-continuous a.s., we can for a.a. $\omega \in \Omega$ define its derivatives with respect to $t$ and $x$ in the sense of distributions.
 Thus we define the time-space white noise $\overset{\bullet}{B}(t,x)=W(t,x)=W(t,x,\omega)$ by
 
\begin{equation}
  	\overset{\bullet}{B}(t,x)=W(t,x)=\frac{\partial}{\partial t}\frac{\partial^{d}B(t,x)}{\partial x_{1}...\partial x_{d}}. \label{WN}
  \end{equation}
In particular, for $d=1$ and $x_1=t$ and get the familiar identity
$$\overset{\bullet}{B}(t)=W(t)=\frac{d}{dt}B(t)\hbox{ in } \S^{'}(\mathbb{R}).$$ 

The definition \eqref{WN} can also be interpreted as an element of the Hida space $(\S)^*$ of \emph{stochastic distributions}, and in that setting it has been proved (see Lindstr\o m, \O ksendal \& Ub\o e \cite{LOU} and Benth \cite{B}) that the Ito-Skorohod integral with respect to $B(dt,dx)$ can be expressed as 
\begin{align}
\int_0^T \int_{\mathbb{R}^d} f(t,x,\omega) B(dt, dx)=\int_0^T \int_{\mathbb{R}^d} f(t,x,\omega) \diamond W(t, x) dt dx,
\end{align}
where $\diamond$ denotes the Wick product.\\
\emph{This is the interpretation we are using in \eqref{heat1} and elsewhere in this paper.}

In particular, if $f(t, x,\omega)=f(t, x)$ is deterministic, we can drop the Wick product and get
\begin{align}
\int_0^T \int_{\mathbb{R}^d} f(t, x) B(dt, dx)=\int_0^T \int_{\mathbb{R}^d} f(t,x)  W(t,x) dt dx.
\end{align}

\section{Fractional Brownian motion and fractional white noise calculus}

There are several ways to introduce fractional Brownian motion and fractional white noise calculus. One commonly used approach due to Hu et al \cite{Hu} is recalled in the Appendix. For our purpose it is more convenient to follow the approach presented by, for example,  Elliott and van der Hoek (see \cite{EvdH} and the references therein, including \cite{L}). Hence we introduce the multi-parameter (time-space) fractional Brownian motion  as obtained by applying an operator $M$ to the classical multi-parameter Brownian motion.  This approach has several advantages. For example, all the different fBm's corresponding to different Hurst coefficients $H \in (0,1)$ are defined on the same white noise probability space $(\Omega, \mathbb{F},\P)$ introduced in Section 3. There are also several computational advantages, as shown later. See also e.g. \cite{DO}.\\

In the following we let $f: \mathbb{R}^n \mapsto \mathbb{R}$ be a function such that its Fourier transform
\begin{equation}
\mathcal{F}(f)(y):= \widehat{f}(y):= (2 \pi)^{-\tfrac{n}{2}} \int_{\mathbb{R}^n} e^{-ixy}f(x)dx; \quad y\in \mathbb{R}^n,\quad ( i=\sqrt{-1})
 \end{equation}
exists.
Let $\mathcal{F}^{-1}$ denote the inverse Fourier transform, defined by
 \begin{align}
 \mathcal{F}^{-1} g(x):= (2 \pi)^{-\tfrac{n}{2}}\int_{\mathbb{R}^n} e^{ixy} g(y) dy,
 \end{align}
 for all functions $g$ such that the integral converges.

\begin{definition}(\cite{EvdH},(A10))
Let $H=(H_1, H_2, ... ,H_n) \in (\tfrac{1}{2},1)^n $ and let $\mathcal{S}=\mathcal{S}(\mathbb{R}^n)$ be the Schwartz space of rapidly decreasing smooth functions on $\mathbb{R}^n$. The operator $M=M^{(H)}: \mathcal{S} \mapsto \mathcal{S}$ is defined by
\begin{align}
Mf(x)&=\mathcal{F}^{-1}(\prod_{j=1}^n |y_j |^{\frac{1}{2}-H_j} (\mathcal{F}f)(y))(x)\\
&=(2 \pi)^{-\tfrac{d}{2}}\int_{\mathbb{R}^n} e^{ixy} (\prod_{j=1}^n |y_j|^{\frac{1}{2} - H_j}(2 \pi)^{-\tfrac{d}{2}} \int_{\mathbb{R}^n} e^{-izy}f(z)dz) dy
\end{align}
If $H_j > \frac{1}{2} $ for all $j$ then the operator $M$ can also be defined as follows:
\begin{align}
M_Hf(x)=  \int_{\mathbb{R}^n} \prod_{j=1}^n C_{H_j}| y_j | ^{H_j-\frac{3}{2}} f(x+y) dy; \quad f \in \mathcal{S},\label{MHa}
\end{align}
where
\begin{equation}
C_{H_j}=[2\Gamma(H_j-\tfrac{1}{2}) \cos(\tfrac{\pi}{2}(H_j-\tfrac{1}{2}))] ^{-1},
\end{equation}
$\Gamma$ being the classical Gamma function.
\end{definition}
\emph{From now on we assume that $H \in (\frac{1}{2},1)^n$.}\\

Since $H \in (\frac{1}{2},1)^n$ is fixed, we will suppress the index $H$ and write $M^{(H)} = M$ in the following.\\
We can in a natural way extend the operator $M$  from $\mathcal{S}$ to the space
\begin{align}
L_{H}^2(\mathbb{R}^n):= \{ f: \mathbb{R}^n \mapsto \mathbb{R}; Mf \in L^2(\mathbb{R}^n )\}.
\end{align}
Note that $L_{H}^2 (\mathbb{R}^n)$ is a Hilbert space when equipped with the inner product
\begin{equation}
(f,g)_H := (f,g)_{L_H^2(\mathbb{R}^n)}=(Mf,Mg)_{L^2 (\mathbb{R}^n)};\quad f,g \in L^2_{H}(\mathbb{R}^n).
\end{equation}

The following result is useful:
\begin{lemma}
Let $\varphi \in L_H^2(\mathbb{R}^n)$ and $\psi \in L^2_{H}(\mathbb{R}^n), \varphi \geq 0, \psi \geq 0$. Then
\begin{align}
(\varphi,M\psi)_{L^2(\mathbb{R}^n)}=(M\varphi,\psi)_{L^2(\mathbb{R}^n)}
\end{align}

\end{lemma}

\dproof
Since
\begin{equation}
Mf(x)= \int_{\mathbb{R}^n} \prod_{j=1}^n C_{H_j}| y_j | ^{H_j-\frac{3}{2}} f(x+y) dy; \quad f \in \mathcal{S},\label{MH}
\end{equation}
we get by the Fubini-Tonelli theorem
\begin{align*}
(\varphi, M\psi)_{L^2(\mathbb{R}^n)}&= \int_{\mathbb{R}^n} \varphi(x)M\psi(x)dx= \int_{\mathbb{R}^n}\varphi(x) \int_{\mathbb{R}^n} \prod_{j=1}^n C_{H_j}| y_j | ^{H_j-\frac{3}{2}} \psi(x+y) dy dx\\
&= \int_{\mathbb{R}^n}\int_{\mathbb{R}^n}\varphi(z-y)\prod_{j=1}^n C_{H_j}| y_j | ^{H_j-\frac{3}{2}}\psi(z) dy dz\\
&=  \int_{\mathbb{R}^n}\int_{\mathbb{R}^n}\psi(z) \prod_{j=1}^n C_{H_j}| y_j | ^{H_j-\frac{3}{2}} \varphi(z+y) dy dz=(M\varphi,\psi)_{L^2(\mathbb{R})},
\end{align*}
which completes the proof.

\fproof

The following property of the operator $M$ will be used in the sequel.

\begin{lemma}\label{l1} 
The operator $M$ satisfies following property:
\begin{align} \label{eq3.5a}
M^2f(x)= \int_{\mathbb{R}^n} \int_{\mathbb{R}^n} \prod_{j=1}^n C_{H_j}^2| y_j z_j| ^{H_j-\frac{3}{2}}  f(y+z+x)dy dz;\quad f \geq 0.
\end{align}
\end{lemma}

\dproof 
Since
$$Mf(x)=\int_{\mathbb{R}^n}\prod_{j=1}^n C_{H_j}| y_j | ^{H_j-\frac{3}{2}} f(x+y)dy$$
we get
\begin{align*}
M^2f(x)&=M(Mf)(x)= \int_{\mathbb{R}^n}\prod_{j=1}^n C_{H_j}| z_j | ^{H_j-\frac{3}{2}}  Mf(x+z) dz \\
&= \int_{\mathbb{R}^n} \prod_{j=1}^n C_{H_j}| z_j | ^{H_j-\frac{3}{2}} \int_{\mathbb{R}^n}\prod_{j=1}^n C_{H_j}| y_j | ^{H_j-\frac{3}{2}}f(x+z+y)dy dz\\
&=\int_{\mathbb{R}^n} \int_{\mathbb{R}^n}\prod_{j=1}^n C_{H_j}^2| y_j z_j | ^{H_j-\frac{3}{2}} f(y+z+x)dy dz.
\end{align*}
This completes the proof.
\fproof

Results from previous lemmas can be gathered in the extended equality which is given in following lemma.

\begin{lemma} \label{5.4} Let  $f,g:\mathbb{R}^n\rightarrow \mathbb{R}_{+}$. Then
	\begin{align} \label{eq l21}
	\int_{\mathbb{R}^n}f(x)	M^2 g(x)dx= \int_{\mathbb{R}^n}g(x)M^2 g(x)dx=\int_{\mathbb{R}^n} Mf(x) Mg(x) dx.
	\end{align}
\end{lemma}

\dproof 
From Lemma \ref{l1} it follows that
\begin{align}
\int_{\mathbb{R}^n} f(x)	M^2 g(x)dx&=\int_{\mathbb{R}^n}f(x) \int_{\mathbb{R}^n} \int_{\mathbb{R}^n} \prod_{j=1}^n C_{H_j}^2| y_j z_j| ^{H_j-\frac{3}{2}}  g(y+z+x)dy dz dx\label{eq3.5}
\end{align}.

If we substitute

\begin{equation*} 
	\left\{ 
        \begin{array}{l} 
            u=x+y+z\\
            v=y\\
            w=z\\
        \end{array}
        \right. 	\Leftrightarrow  \begin{array}{l} 
            x=u-v-w=: X(u,v,w)\\
            y=v=:Y(u,v,w)\\
            z=w=:Z(u,v,w),\\
        \end{array}
    \end{equation*}    
    then the Jacobian is 
    \begin{equation}
    \begin{vmatrix} 
\frac{\partial X}{\partial u} & \frac{\partial X}{\partial v} & \frac{\partial X}{\partial w} \\ 
\frac{\partial Y}{\partial u} & \frac{\partial Y}{\partial v} & \frac{\partial Y}{\partial w} \\ 
\frac{\partial Z}{\partial u} & \frac{\partial Z}{\partial v} & \frac{\partial Z}{\partial w} 
\end{vmatrix} =   
    \begin{vmatrix}
1 & -1 & -1 \\ 
0 & 1 & 0\\ 
0 & 0 & 1
\end{vmatrix}=1,
\end{equation}	

and we obtain
\begin{align}
&\int_{\mathbb{R}^n}f(x)	M^2g(x) dx=\int_{\mathbb{R}^n}\int_{\mathbb{R}^n}\int_{\mathbb{R}^n}\prod_{j=1}^n C_{H_j}^2| u_j w_j| ^{H_j-\frac{3}{2}}  f(u-v-w)g(u)du dv dw\nonumber\\
&=\int_{\mathbb{R}^n}\int_{\mathbb{R}^n}\int_{\mathbb{R}^n} g(u)\prod_{j=1}^n C_{H_j}^2| u_j w_j| ^{H_j-\frac{3}{2}}  f(u-v-w)g(u)du dv dw\nonumber\\	&=\int_{\mathbb{R}^n}g(u)	M^2 f(u)du.
\end{align}
This proves the first identity. A similar argument proves the second identity.
\fproof

\medskip

Now we can give a formal definition of multi-parameter fractional Brownian motion:

\begin{definition}{(Fractional Brownian motion revisited)}\\
For $x \in \mathbb{R}^n$ define
\begin{equation} \label{eq4.5a}
 B_{H}(x):= B_{H}(x,\omega):=\langle \omega, M_x(\cdot) \rangle ; \quad \omega \in \Omega 
\end{equation}
where we have put 
\begin{align}
M_x (y):=M(\chi_{[0,x_1]\times [0,x_2] \times ... \times [0,x_n]})(y)=\prod_{j=1}^n C_{H_j}^{-1}\Big(\frac{x_j-y_j}{|x_j-y_j|^{\frac{3}{2}-H}}+\frac{y_j}{|y_j|^{\frac{3}{2}-H}}\Big);\quad y,x \in \mathbb{R}^n.
\end{align}

Then $B_{H}$ is a Gaussian process, and by the properties of the classical Bm we have 
\begin{equation}
\E[B_{H} (x)]=B_{H}(0)=0,
\end{equation}
and (see \cite{EvdH},(A10))
\begin{equation}
\E[B_{H}(x) B_{H}(y)]=  \prod_{j=1}^{d}C_{H_j}\big( |x_{j}|^{2H_{j}}+|y_{j}|^{2H_{j}}-|x_{j}-y_{j}|^{2H_{j}}\big); \quad x,y \in \mathbb{R}^n,
\end{equation}
where $\E$ denotes expectation with respect to $\P$.\\
By the Kolmogorov continuity theorem the process $B_{H}$ has a continuous version, which we will also denote by $B_{H}.$ Except for a constant, this process coincides with the fractional Brownian motion $B_H(x)$ introduced in the Appendix.
From now on we call this process 
 \emph{ the fractional Brownian motion (fBm) with Hurst parameter $H$} and denote it by $B_H(x)$.
\end{definition}

\subsection{Fractional white noise calculus}
By \eqref{eq4.5a} it follows by approximating $f \in \mathcal{S}(\mathbb{R}^n)$ by step functions that
\begin{equation}\label{int}
 \int_{\mathbb{R}^n} f(x) dB_{H}(x)=\langle \omega,Mf \rangle =\int_{\mathbb{R}^n} Mf(x) dB(x)
\end{equation}
and hence
\begin{equation}
B_{H}(x)=\int_{\mathbb{R}^n} \prod_{j=1}^n \chi_{[0,x_j]}(y_j) dB_{H}(y)= \int_{\mathbb{R}^n}M_{x}(y)dB(y)=\langle \omega, M_x \rangle.
\end{equation}
We know that the action of $\omega \in \mathcal{S}'(\mathbb{R}^n)$ extends from $\mathcal{S}(\mathbb{R}^n)$ to $L^2(\mathbb{R}^n),$ and that if $f \in \mathcal{S}$ then $Mf \in L^2(\mathbb{R}^n).$ Therefore we can define $M$ on $\mathcal{S}'(\mathbb{R}^n)$ by setting
\begin{equation}
\langle M\omega ,\varphi \rangle= \langle \omega, M\varphi \rangle; \varphi \in \mathcal{S}(\mathbb{R}^n), \omega \in \mathcal{S}'(\mathbb{R}^n).
\end{equation}
Using this, we can write
\begin{align}
B_{H}(x,\omega)=\langle \omega, M_x \rangle= \langle M\omega, \prod_{j=1}^n\chi_{[0.x_j](\cdot)} \rangle =B(x,M\omega).
\end{align}

We define the \emph{fractional white noise} 
$W_{H}(x)=W_{H}(x,\omega)$ by
\begin{align}
W_{H}(x,\omega)= \frac{\partial^{n}}{\partial x_{1}...\partial x_{n}}B_{H}(x) \text{ (derivative in } (\mathcal{S})^{*} ).
\end{align}
By the above it follows that 

\begin{align}
W_{H}(x,\omega)=W(x,M\omega),
\end{align}
 where 
\begin{align}
W(x)=W(x,\omega)= \frac{\partial^{n}}{\partial x_{1}...\partial x_{n}}B(x)  \text{ is the classical white noise in } (\mathcal{S})^{*}.
\end{align}
In general we define the \emph{Wick-Ito-Skorohod} (WIS)-integral  as follows:

\begin{definition}  \label{wis}
(The fractional Wick--It\^o--Skorohod (WIS) integral) \\
Let  $f : \mathbb{R}^n \rightarrow (\mathcal{S})^{*}$ be such that $f (x)  \diamond W_{H}(x) $ is $dx$-integrable in $ (\mathcal{S})^{*}$. Then we say that  $f$ is
$dB_{H}$-integrable and we define the Wick--It\^o--Skorohod (WIS) integral of $f (\cdot) = f(\cdot, \omega)$ with respect to $B_{H}$ by
\begin{align}
\int_{\mathbb{R}^n} f(x) dB_{H}(x)= \int_{\mathbb{R}^n} f(x) \diamond W_{H}(x) dx,
\end{align}
where the integral on the right is the (Pettis) Lebesgue integral with values in $(\mathcal{S})^{*}.$
\end{definition}

By \eqref{int} and Lemma \ref{5.4} we obtain the following isometry:
\begin{theorem} \label{iso1}
Assume that $f \geq 0$ is deterministic. Then

\begin{align}
    \E[(\int_{\mathbb{R}^n} f(x) dB_{H}(x))^2]= \E[\int_{\mathbb{R}^n} (Mf)^2(x) dx]= \E[\int_{\mathbb{R}^n} f(x) M^2f(x) dx].
    \end{align}
\end{theorem}

\section{The fractional heat equation with additive noise}

In this section  we study the additive noise case, i.e. we study the equation  
\begin{align}
 	 &\frac{\partial^{\alpha}}{\partial t^{\alpha}}Y(t,x)=-\lambda (-\Delta)^{ \beta} Y(t,x)+\langle b,  \nabla Y(t,x)\rangle+\sigma W_{H}(t,x)\; (t,x)\in (0,\infty)\times \mathbb{R}^{d} \label{frac}
	 \end{align}
   with boundary conditions
\begin{align}\label{add01}
    Y(0,x)&=\delta_0(x)\text{ (the Dirac measure at  } x), \nonumber\\
    \lim_{x \rightarrow +/- \text{ }\infty}Y(t,x)&=0,
\end{align}
  \begin{definition}
   We say that a process $Y(t,x)=Y(t,x,\omega): [0,T] \times \mathbb{R}^d \times \Omega \mapsto \mathbb{R}$ is a solution of of the fractional stochastic heat equation \eqref{frac} - \eqref{add01} if \\
(i) $\omega \mapsto Y(t,x,\omega) \in (\mathcal{S})^{*}$ for all $t,x$, \\
(ii) $(t,x) \mapsto Y(t,x,\omega) \in \mathcal{S}'(\mathbb{R}^{1+d})$ for all $\omega \in \Omega$, and\\
(iii) the function $(t,x) \mapsto Y(t,x,\omega)$ satisfies the equation \eqref{frac} - \eqref{add01} in the sense of distribution, for all $\omega$.
\end{definition}  
Our first main result in this section is the following:
 \begin{theorem} \label{add}
     The unique solution $Y(t,x) \in \mathcal{S}'$ of the fractional heat equation \eqref{frac} - \eqref{add01} is given by
     \begin{align}\label{addsol}
     &Y(t,x)=(2\pi)^{-d} \int_{\mathbb{R}^d} e^{ixy}E_{\alpha}([-\lambda |y|^{2\beta}+i\langle b,y \rangle] t^{\alpha})dy \nonumber\\
    &+ \sigma (2\pi)^{-d} \int_{0}^{t}(t-r)^{\alpha -1}\int_{\R^d} \Big(\int_{\mathbb{R}^{d}}e^{i(x-z)y}E_{\alpha,\alpha}([-\lambda |y|^{2\beta}+i\langle b,y \rangle] (t-r)^{\alpha})dy\Big) B_H(dr,dz)
    \end{align}
    \textcolor{blue}{For $\alpha=2$, this formula assumes the initial velocity $\partial_tY(0,x)=0$; for $\alpha<2$, only the initial value is needed.}
 \end{theorem}  
 \dproof
 In the following we let
\begin{align}
    L\{f(\cdot,x)\}(s,x)=\widetilde{f(\cdot,x)}(s,x)=\int_0^t  e^{-rs}f(r,x)dr;\quad s\geq 0
\end{align}
denote the Laplace transform of $f(t,x)$ as a function of $t$, and we let

\begin{align}
    F\{f(t,\cdot)\}(t,y)=\widehat{f(t,\cdot)}(t,y)
    \end{align}
    denote the (multi-dimensional) Fourier transform of $f(t,x)$ as a function of $x=(x_1, x_2, ... ,x_d).$ \\

First assume that $Y(t,x)$ is a solution of \eqref{frac}--\eqref{add01}. We apply the Laplace transform $L$
to both sides of  \eqref{frac} and obtain

\begin{equation}\label{step1b}
	s^{\alpha}\widetilde{Y}(s,x)-s^{\alpha-1}Y(0^{+},x)=-\widetilde{\lambda (-\Delta)^{ \beta} Y}(s,x)+\widetilde{\langle b, \nabla Y\rangle}(s,x)+\sigma \widetilde{W}_H(s,x).
	\end{equation}
Applying the Fourier transform $F$ to this equation  	 we get, since $\widehat{Y}(0^{+},y)=1$,
\begin{equation}\label{step1a}
s^{\alpha}\widehat{\widetilde{Y}}(s,y)-s^{\alpha-1}=-\widehat{\widetilde{\lambda (-\Delta)^{ \beta} Y}}(s,y)+\widehat{\widetilde{\langle b,\nabla Y\rangle}}(s,y)+\sigma \widehat{\widetilde{W_H}}(s,x).
	\end{equation}
Since the Laplace and Fourier transforms commute, we get
\begin{equation}\label{step1a2}
s^{\alpha}\widehat{\widetilde{Y}}(s,y)-s^{\alpha-1}=-\widetilde{\widehat{\lambda (-\Delta)^{ \beta} Y}}(s,y)+\widetilde{\widehat{\langle b,\nabla Y\rangle}}(s,y)+\sigma \widetilde{\widehat{W_H}}(s,x).
	\end{equation}
   
   Using that $$\widehat {(-\Delta)^{\beta}Y}(s,y)= |y|^{2\beta}$$ and $$\widehat{\nabla Y}(s,y)=i y \widehat{Y}(s,y),$$ we obtain
   
 \begin{equation}\label{step1a3}
s^{\alpha}\widetilde{\widehat{Y}}(s,y)-s^{\alpha-1}=-\lambda \widetilde{|y|^{2\beta}\widehat{Y}}(s,y)+i \widetilde{\langle b,y\rangle\widehat{Y}}(s,y)+\sigma \widetilde{\widehat{W_H}}(s,x).
	\end{equation}  
   Hence
\begin{equation}\label{step2a}
   	\left( s^{\alpha}+\l|y|^2 - i\langle b,y\rangle \right)\widetilde{\widehat{Y}}(s,y)=  s^{\alpha-1}+\sigma\widetilde{\widehat{W_H}}(s,y)
   \end{equation}
or
   \begin{align}\label{step3}
       \widehat{\widetilde{Y}}(s,y)=\frac{s^{\alpha -1}}{s^{\alpha} + \l |y|^{2\beta} -i\langle b,y\rangle}
       + \frac{\sigma \widehat{\widetilde{W}_H}(s,y)}{s^{\alpha} + \l |y|^{2\beta} -i\langle b,y\rangle}.
   \end{align}
   
 Applying the inverse Laplace operator $L^{-1}$  to this equation we get, using \eqref{L2},
 \begin{align} 
       \widehat{Y}(t,y)&=L^{-1} \Big(\frac{s^{\alpha -1}}{s^{\alpha} + \l |y|^{2\beta} -i\langle b,y\rangle}\Big)(t,y)
       + L^{-1}\Big(\frac{\sigma \widetilde{\widehat{W}_H}(s,y)}{s^{\alpha} + \l |y|^{2\beta} -i\langle b,y\rangle}\Big)(t,y)\nonumber\\
       &=E_{\alpha,1}((-\l |y|^{2\beta} + i\langle b,y \rangle) t^{\alpha}) + L^{-1}\Big(\frac{\sigma \widetilde{\widehat{W}_H}(s,y)}{s^{\alpha} + \l |y|^{2\beta} -i\langle b,y\rangle}\Big)(t,y),\label{step 5}
   \end{align}
   where in general
   \begin{align} \label{add2.12}
   E_{\alpha,\beta}(z)=\sum_{k=0}^{\infty} \frac{z^{k}}{\Gamma(\alpha k + \beta)}    
   \end{align}
   is the Mittag-Leffler function.  \vskip 0.2cm
   It remains to find 
$L^{-1}\left(\frac{\sigma\widehat{\widetilde{W}_H}(s,y)}{s^{\alpha} + \l |y|^{2\beta} -i\langle b,y\rangle}\right)$:\\
Recall that the convolution $f\ast g$ of two functions $f,g: [0,\infty) \mapsto \mathbb{R}$  is defined by
\begin{align}
(f \ast g)(t)=\int_0^t f(t-r)g(r) dr; \quad t \geq 0.
\end{align}
The convolution rule for Laplace transform states that $$L\left( \int_{0}^{t}f(t-r)g(r)dr\right) (s)=Lf(s)Lg(s),$$ 
or 
\begin{equation}\label{addititive 12}
	\int_{0}^{t}f(t-w)g(w)dw=L^{-1}\left( Lf(s)Lg(s)\right) (t).\\
\end{equation}
Using \eqref{L3} we get
\begin{align*}
	&L^{-1}\left(\frac{1}{s^{\alpha} + \l |y|^{2\beta} -i\langle b,y\rangle} \right) (t)\\
    &=t^{\alpha-1}E_{\alpha,\alpha}([\l t^{\alpha}|y|^{2\beta} -i \langle b,y \rangle] t^{\alpha})
	=\sum_{k=0}^{\infty}\frac{t^{\alpha-1}([-\l t^{\alpha}|y|^{2\beta}+i\langle b,y \rangle] t^{\alpha})^{k}}{\Gamma(\alpha k+\alpha)}\nonumber\\
	&=\sum_{k=0}^{\infty}\frac{(-\l t^{\alpha}|y|^{2\beta}+i\langle b,y \rangle)^{k}t^{\alpha(k+1)-1}}{\Gamma(\alpha(k+1))}
	=: \Lambda(t,y).\label{addL}
\end{align*}
In other words,
\begin{equation}
	\frac{\sigma}{s^{\alpha} + \l |y|^{2\beta} -i\langle b,y\rangle}=\sigma L \Lambda(s,y).
\end{equation}
Combining this with \eqref{addititive 12} we get 
\begin{align}
	L^{-1}\left( \frac{\sigma}{s^{\alpha} + \l |y|^{2\beta} -i\langle b,y\rangle} \widehat{\widetilde{W}} (s,y)\right) (t)&=L^{-1}\left( L\left( \sigma\Lambda(s,y)\right) \widetilde{\widehat{W}}(s,y)\right) (t)\\
	&=\sigma \int_{0}^{t}\Lambda(t-r,y)\widehat{W}_H(r,y)dr.
\end{align}
Substituting this into \eqref{step 5} we get 
\begin{equation}
	\widehat{Y}(t,y)=E_{\alpha,1}\left( [-\l |y|^{2\beta} +i \langle b,y \rangle] t^{\alpha}\right)+\sigma\int_{0}^{t}\Lambda(t-r,y)\widehat{W}_H(r,y)dr.
\end{equation}
Taking inverse Fourier  transform we end up with
\begin{equation}\label{add2.18a}
 Y(t,x)=F^{-1}\left( E_{\alpha,1}\left([-\l |y|^{2\beta} +i\langle b,y \rangle]t^{\alpha} 
 \right)\right)(x)+\sigma F^{-1}\left(\int_{0}^{t}\Lambda(t-r,y)\widehat{W}_H(r,y)dr\right)(x).
\end{equation}
 
Now, we use that 
\begin{align}
F\{\int_{\R^d}f(x-z)g(z)dz\}(y)= F\{f\}(y)F\{g\}(y),\quad x \in \mathbb{R}^d
\end{align}
or
\begin{equation}
	\int_{\R^d}f(x-z)g(z)dz=F^{-1}\Big( F\{f\}(y)F\{g\}(y)\Big)(x)\quad x \in \mathbb{R}^d,
\end{equation}
 where $x=(x_1, x_2, ... , x_d), z=(z_1, z_2, ... ,z_d)$,

This gives,
\begin{align*}
	&F^{-1}\left( \int_{0}^{t}\Lambda(t-r,y)\widehat{W}(r,y)dr\right) (x)=\int_{0}^{t}F^{-1}\left( \Lambda(t-r,y)\widehat{W}_H(r,y)\right) (x)dr\\
	&=\int_{0}^{t}F^{-1}\Big( F\left( F^{-1}\Lambda(t-r,y)\right) (y)F (W_H(r,x))(y)\Big) (x)dr\\
	&=\int^t _{0}\int_{\R^d}\left( F^{-1}\Lambda(t-r,y)(x-z)\right) W_H(r,z)dzdr\\
	&=\int_{0}^{t}\int_{\R^d}\left( (2\pi)^{-d}\int_{\mathbb{R}^{d}}e^{i(x-z)y}\Lambda(t-r,y)dy\right) W_H(r,z)dzdr\\
	&=(2\pi)^{-d}\int_{0}^{t}\int_{\R^d} \left( \int_{\mathbb{R}^{d}}e^{i(x-z)y}\Lambda(t-r,y)dy\right) B_H(dr,dz).
\end{align*}

Combining this with \eqref{add2.18a}, \eqref{add2.12}, and the definition of
$\Lambda$ above, we get
\begin{align}\label{addexp1a}
    &Y(t,x)= F^{-1}(\sum_{k=0}^{\infty} \frac{([- \l |y|^{2\beta} + i \langle b,y \rangle]t^{\alpha})^k}{\Gamma(\alpha k +1)})+\sigma (2\pi)^{-d} \int_{0}^{t}\int_{\R^d} \left(\int_{\mathbb{R}^{d}}e^{i(x-z)y}\Lambda(t-r,y)dy\right) B_H(dr,dz)\nonumber\\
    &=(2\pi)^{-d} \int_{\mathbb{R}^d} e^{ixy}\sum_{k=0}^{\infty} \frac{([- \l |y|^{2\beta} + i \langle b,y \rangle]t^{\alpha})^k}{\Gamma(\alpha k +1)}dy \nonumber\\
    &+ \sigma (2\pi)^{-d} \int_{0}^{t}(t-r)^{\alpha -1}\int_{\R^d} \left(\int_{\mathbb{R}^{d}}e^{i(x-z)y}\sum_{k=0}^{\infty}\frac{(-\l |y|^{2\beta}+i\langle b,y \rangle)^{k}(t-r)^{\alpha(k+1)-1}}{\Gamma(\alpha(k+1))}dy\right) B_H(dr,dz) \\
    &=(2\pi)^{-d} \int_{\mathbb{R}^d} e^{ixy}E_{\alpha}([-\lambda |y|^{2\beta}+i\langle b,y \rangle] t^{\alpha})dy \nonumber\\
    &+ \sigma (2\pi)^{-d} \int_{0}^{t}(t-r)^{\alpha -1}\int_{\R^d} \Big(\int_{\mathbb{R}^{d}}e^{i(x-z)y}E_{\alpha,\alpha}([-\lambda |y|^{2\beta}+i\langle b,y \rangle] (t-r)^{\alpha})dy\Big) B_H(dr,dz)
\end{align}
This proves uniqueness and also that the unique solution (if it exists) is given by \eqref{addsol}.\\

b) Next, define $Y(t,x)$ by \eqref{addsol}. Then we see that we can prove that $Y(t,x)$ satisfies \eqref{add} by reversing the argument above.
We skip the details.
\fproof

\subsection{The classical case ($\alpha$ = 1)}
It is interesting to compare the above result with the classical case when $\alpha$=1:\\
If $\alpha=1$, we get $Y(t,x)=I_{0}(x)+I_{1}(t,x)$, where
\begin{align*}
	&I_{0}(x)=(2\pi)^{-d}\int_{\mathbb{R}^{d}}e^{ixy}\sum_{k=0}^{\infty}\frac{\left( -\l t |y|^{2}\right) ^{k}}{k!}dy=(2\pi)^{-d} \int_{\mathbb{R}^d} e^{ixy}E_{\alpha}(-\lambda |y|^2 t^{\alpha})dy, \\
&I_1(t,x)= \sigma (2\pi)^{-d} \int_{0}^{t}\int_{0}^x \Big(\int_{\mathbb{R}^{d}}e^{i(x-z)y}E_{\alpha}(-\lambda |y|^2 (t-r)^{\alpha})\Big) B_H(dr,dz)
\end{align*}
Using that that $\Gamma(k+1)=k!$ 
and the Taylor expansion of the exponential function, we get
\begin{align*}
 I_{0}=(2\pi)^{-d}\int_{\mathbb{R}^{d}}e^{ixy}e^{-\l t |y|^{2}} dy
 =(2\pi)^{-d}\left( \frac{\pi}{\l t}\right) ^{\frac{d}{2}}e^{-\frac{|x|^{2}}{4\l t}}
 =(4 \pi \l t) ^{-\frac{d}{2}} e^{-\frac{|x|^{2}}{4\l t}},	
\end{align*}
where we used the general formula
\begin{equation}
\int_{\mathbb{R}^{d}}e^{-\left( a |y|^{2}+2by\right) }dy=\left(\frac{\pi}{a} \right)^{\frac{d}{2}} e^{\frac{b^{2}}{a}};\;a>0;\;b\in\mathbb{C}^d.\label{exp}
\end{equation}\\
Similarly,
\begin{align*}
	I_{1}(t,x)&=\sigma(2\pi)^{-d}\int_{0}^{t}\int_{0} ^x \int_{\mathbb{R}^{d}}e^{i(x-z)y}\sum_{k=0}^{\infty}\frac{\left(-\l (t-r)|y|^{2} \right)^{k} }{k!}dyB_H(dr,dz)\\
	&=\sigma(2\pi)^{-d}\int_{0}^{t}\int_{0}^x \left( \frac{\pi}{\l (t-r)}\right) ^{\frac{d}{2}}e^{-\frac{|x-z|^{2}}{4\l (t-r) }}B_H(dr,dz)\\
	&=\sigma (4 \pi \l)^{-\tfrac{d}{2}}\int_0^t \int_{0}^x  (t-r)^{-\tfrac{d}{2}} e^{-\frac{|x-z|^2}{4 \l (t-r)}} B_H(dr,dz).
\end{align*}
Summarizing the above, we get, for $\alpha=1$, 
\begin{align}
	Y(t,x)&=(4\pi \l t)^{-\frac{d}{2}}e^{-\frac{|x|^{2}}{4\l t}}\nonumber\\
	&+\sigma(4\pi\l)^{-\frac{d}{2}}\int_{0}^{t}\int_{0}^x (t-r)^{-\frac{d}{2}}e^{-\frac{|x-z|^{2}}{ 4\l (t-r)}}B_H(dr,dz)\;\; \label{alpha1}
\end{align}
This is in agreement with a well-known classical result. See e.g. Section 4.1 in \cite{Hu}.

\subsection{When is $Y(t,x)$ a mild solution?}

It was pointed out already in 1984 by John Walsh \cite{W} that (classical) SPDEs driven by time-space white noise $W(t,x); (t,x) \in [0,\infty) \times \R^d$ may have only distribution valued solutions if $d \geq 2$. 
Indeed, the solution $Y(t,x)$ that we found in the previous section is in general distribution valued. But in some cases the solution can be represented as an element of $L^2(\P)$. Following Y. Hu \cite{Hu} we make the following definition:

\begin{definition}
The solution $Y(t,x)$ is called \emph{mild} if $Y(t,x) \in L^2(\P)$ \\for all $t>0, x >0$.
\end{definition}

We now ask the following:
\begin{problem}
   For what values of $\alpha \in (0,2)$ and what dimensions $d=1,2, ...$ is $Y(t,x)$ mild?
      \end{problem}

\textcolor{blue}{%
\textbf{Correction:} In the original version, the isometry used in the
computation of $J_2$ was the classical It\^o isometry, which does not
account for the fractional nature of the noise. The correct isometry
is the fractional isometry given in Theorem \ref{iso1}.
This correction introduces the Hurst parameter $H_0$ in the mild
condition.
}

To answer this we first note that by the \textcolor{blue}{fractional isometry (Theorem \ref{iso1})} we have 
\begin{equation}
    \E\Big[ Y^{2}(t,x)\Big]=J_{1}+J_{2},
\end{equation}
where
\small
 \begin{align}
 J_1&=(2\pi)^{-2d} \Big(\int_{\mathbb{R}^d} e^{ixy}\sum_{k=0}^{\infty} \frac{(- \l t^{\alpha} |y|^2)^k}{\Gamma(\alpha k +1)}dy\Big)^2 \nonumber\\
   J_2 &=\sigma^2\int_0^t\int_{\R^d} |M K_\alpha(t-r,x-z)|^2\,dr\,dz,
 \end{align}
with $K_\alpha$ defined by its Fourier transform
\begin{equation}
 \widehat K_\alpha(s,\xi)
 =
 s^{\alpha-1}E_{\alpha,\alpha}(-\lambda s^\alpha|\xi|^{2\beta}).
\end{equation}

\textcolor{blue}{%
Using the Fourier multiplier property of $M$:
\[
\widehat{M K_\alpha}(s,\xi)
=
\prod_{j=0}^d |\xi_j|^{1/2-H_j}\,\widehat K_\alpha(s,\xi),
\]
we obtain, after the change of variables $\eta=s^{\alpha/(2\beta)}\xi$,
\begin{equation}
\int_{\R^d} |M K_\alpha(s,\xi)|^2\,d\xi
=
C_{\alpha,\beta,d,\lambda,H_0}
s^{2\alpha-2-\alpha d/(2\beta)-\alpha(1-2H_0)/(2\beta)}.
\end{equation}
The time integral near zero is finite if and only if
\begin{equation}
d+(1-2H_0)<\frac{2\beta(2\alpha-1)}{\alpha}.
\end{equation}
}

This is the same estimate as in Theorem in \cite{MO}. Therefore the proof of that theorem applies here, and we can conclude the following:
 
  \begin{theorem} \label{mild1}
 Let $Y(t,x)$ be the solution of the $\alpha$-fractional stochastic heat equation with additive noise.
Then the following holds for all $H=(H_1,H_2, ... ,H_d) \in (\tfrac{1}{2},1)^d$:
\begin{itemize}
\item
a) If $\alpha = 1$, then $Y(t,x)$ is mild if and only if $d=1$.
\item
b) If $\alpha > 1$ then $Y(t,x)$ is mild if $d=1$ or $d=2$.
\item
c) If $\alpha < 1$ then $Y(t,x)$ is not mild for any  $d.$
\end{itemize}
 \end{theorem} 

\begin{remark}
\begin{itemize}
\item
See Y. Hu \cite{Hu}, Proposition 4.1 for a generalisation of the above result in the case $\alpha=1$.
\item
In the cases $\alpha > 1, d \geq 3$ we do not know if the solution $Y(t,x)$ is mild or not. This is a topic for future research.
\end{itemize}
\end{remark}

\textcolor{blue}{
\subsection{The case $\alpha =2$ (wave equation)}
The results above are in particular valid for $\alpha=2$ and $\beta=1, b=0$, in which case our equation \eqref{frac} becomes a stochastic wave equation of the form
\begin{align}
 	 &\frac{\partial^{2}}{\partial t^{2}}Y(t,x)=\lambda \Delta Y(t,x)+\sigma W_{H}(t,x); \ \ (t,x)\in (0,\infty)\times \mathbb{R}^{d} \label{vawe}
	 \end{align}
   with boundary conditions
\begin{align}\label{add02}
    Y(0,x)&=\delta_0(x)\text{ (the Dirac measure at  } x), \nonumber\\
    \lim_{x \rightarrow +/- \text{ }\infty}Y(t,x)&=0,
\end{align}
and the solution is
\begin{align}\label{addsol0}
     &Y(t,x)=(2\pi)^{-d} \int_{\mathbb{R}^d} e^{ixy}E_{2}(-\lambda |y|^{2} t^{2})dy \nonumber\\
    &+ \sigma (2\pi)^{-d} \int_{0}^{t}(t-r)\int_{\R^d} \Big(\int_{\mathbb{R}^{d}}e^{i(x-z)y}E_{2,2}(-\lambda |y|^{2} (t-r)^{2})dy\Big) B_H(dr,dz).
    \end{align}
    }

    \begin{remark}
        Note that this solution formula was found without any restrictions on the space dimension $d$. This may look surprising, in view of the fact that for the classical wave equation the solution formula depends heavily on the dimension $d$.  
        However, there is no contradiction here. Formula (5.31) is dimension-independent only at the Fourier-transform level; the dimension enters through the inverse Fourier transform. When $\alpha$ = 2 and $c=\sqrt{\lambda}$, the Mittag-Leffler functions become the usual cosine and sine wave propagators (see Section \ref{alpha2} below). Their inverse transforms yield the familiar  d \' Alembert, Poisson, and Kirchhoff kernels in dimensions one, two, and three, respectively. Thus, \eqref{addsol0} is valid in every dimension as a tempered-distribution identity, although its physical-space form depends strongly on the dimension.
    \end{remark}
    
\section{The fractional diffusion--wave equation with multiplicative noise}
The  multiplicative problem requires additive conditions when
$1<\alpha\leq 2$.  Because  in this case, the Caputo derivative contains two initial
traces.  In particular, the endpoint $\alpha=2$ is a stochastic wave equation,
not a heat equation.
\begin{align} 
 	 &\frac{\partial^{\alpha}}{\partial t^{\alpha}}Y(t,x)=-\lambda (-\Delta)^{ \beta} Y(t,x)+b \bullet \nabla Y(t,x)+\sigma Y(t,x)\diamond W(t,x)\; (t,x)\in (0,\infty)\times \mathbb{R}^{d}.
	 \end{align}

Let
\begin{align}
 \Psi(\xi)=\lambda |\xi|^{2\beta}-i\langle b,\xi\rangle,
 \qquad \xi\in\R^d,
\end{align}
and consider equation \eqref{heat1} with initial conditions
\begin{align}\label{mult0}
 Y(0,x)&=u_0(x),\nonumber\\
 \partial_tY(0,x)&=u_1(x),\qquad 1<\alpha\leq2.
\end{align}
We assume that $u_0$ and $u_1$ are
bounded measurable functions.

Define the deterministic kernels by their Fourier transforms:
\begin{align}
 \widehat{G}_{\alpha}(t,\xi)
 &=E_{\alpha,1}\big(-t^\alpha\Psi(\xi)\big),\label{multG}\\
 \widehat{G}_{\alpha}^{(1)}(t,\xi)
 &=tE_{\alpha,2}\big(-t^\alpha\Psi(\xi)\big),\label{multG1}\\
 \widehat{K}_{\alpha}(t,\xi)
 &=t^{\alpha-1}E_{\alpha,\alpha}
   \big(-t^\alpha\Psi(\xi)\big).\label{multK}
\end{align}
The kernel $G_{\alpha}^{(1)}$ is used only for $\alpha>1$.  Put
\begin{align}\label{multdet}
 J_0(t,x)&=(G_\alpha(t,\cdot)*u_0)(x)\nonumber\\
 &\quad+\mathbf{1}_{\{\alpha>1\}}
 (G_{\alpha}^{(1)}(t,\cdot)*u_1)(x).
\end{align}

\begin{definition}\label{multdef}
An $(\mathcal S)^*$-valued process $Y$ is a mild solution of the
multiplicative equation if
\begin{align}\label{vol}
 Y(t,x)=J_0(t,x)+\sigma\int_0^t\int_{\R^d}
 K_\alpha(t-s,x-z)Y(s,z)\diamond W(s,z)\,dz\,ds.
\end{align}
If the integrand is predictable and square integrable, the last term agrees
with the Walsh integral
\begin{align}\label{multwalsh}
 \sigma\int_0^t\int_{\R^d}
 K_\alpha(t-s,x-z)Y(s,z)\,B(ds,dz).
\end{align}
 
\end{definition}

\begin{theorem}\label{mult1}
Suppose that all transforms and convolutions below are well defined.  Then an
$(\mathcal S)^*$-valued process solves equation \eqref{heat1} with the initial
conditions \eqref{mult0} in the sense of distributions if and only if it
satisfies the mild equation \eqref{vol}.
\end{theorem}

\dproof
Taking the Fourier transform in space and the Laplace transform in time gives
\begin{align}
 &(s^\alpha+\Psi(\xi))\widetilde{\widehat Y}(s,\xi)
 =s^{\alpha-1}\widehat u_0(\xi)
 +\mathbf{1}_{\{\alpha>1\}}s^{\alpha-2}\widehat u_1(\xi)\nonumber\\
 &\hspace{5cm}+\sigma\widetilde{\widehat{Y\diamond W}}(s,\xi).
\end{align}
We use the identities
\begin{align}
 L^{-1}\left\{\frac{s^{\alpha-1}}{s^\alpha+a}\right\}(t)
 &=E_{\alpha,1}(-at^\alpha),\nonumber\\
 L^{-1}\left\{\frac{s^{\alpha-2}}{s^\alpha+a}\right\}(t)
 &=tE_{\alpha,2}(-at^\alpha),\nonumber\\
 L^{-1}\left\{\frac{1}{s^\alpha+a}\right\}(t)
 &=t^{\alpha-1}E_{\alpha,\alpha}(-at^\alpha).
\end{align}
The convolution theorem first in time and then in space yields \eqref{vol}.
Conversely, applying the two transforms to \eqref{vol} recovers the transformed
equation and the initial traces.  This proves the equivalence.
\fproof

\subsection{Existence, uniqueness, and chaos expansion}
For $T>0$ define
\begin{align}\label{multq}
 q_T:=\int_0^T\int_{\R^d}|K_\alpha(s,z)|^2\,dz\,ds.
\end{align}

\begin{theorem}\label{vol3}
Assume that $q_T<\infty$ and that $J_0$ is bounded on
$[0,T]\times\R^d$.  Then \eqref{vol} has a unique predictable mild solution
with locally bounded second moments.  It is obtained by Picard iteration and
has the Wiener chaos expansion
\begin{align}\label{multchaos}
 Y(t,x)=J_0(t,x)+\sum_{n=1}^{\infty}\sigma^n I_n
 \big(\widetilde f_n(\,\cdot\,;t,x)\big),
\end{align}
where $\widetilde f_n$ denotes symmetrisation and, on
$0<s_1<\cdots<s_n<t$,
\begin{align}
 f_n(\mathbf s,\mathbf z;t,x)
 &=K_\alpha(t-s_n,x-z_n)
   \prod_{j=1}^{n-1}K_\alpha(s_{j+1}-s_j,z_{j+1}-z_j)\nonumber\\
 &\quad\times J_0(s_1,z_1).
\end{align}
\end{theorem}

\dproof
Starting from $Y_0=J_0$, we obtain $Y_{n+1}$ by inserting $Y_n$ in   in \eqref{multwalsh}.  The Walsh isometry and \eqref{multq} give the
usual Volterra estimates for the differences $Y_{n+1}-Y_n$.  

Since the
$L^2$-mass of $K_\alpha$ on $[0,\delta]\times\R^d$ tends to zero with
$\delta$, the Picard map is a contraction on sufficiently short time
intervals.  Repeating  many times covers $[0,T]$.
The same estimate apply to the difference of two solutions proves
uniqueness.  Expanding the Picard iterates gives \eqref{multchaos}; compare
the general random-field construction in \cite{CHN,W}.
\fproof

\subsection{When the solution is mild?}
The condition is the square integrability of the forcing kernel
$K_\alpha$. 
Plancherel's identity gives
\begin{align}\label{multdalang}
 q_T=(2\pi)^{-d}\int_0^T\int_{\R^d}
 s^{2\alpha-2}\left|E_{\alpha,\alpha}
 \big(-s^\alpha\Psi(\xi)\big)\right|^2\,d\xi\,ds.
\end{align}

\begin{theorem}\label{mild2}
Assume $b=0$, $\lambda>0$, and $0<\beta\leq1$.
\begin{itemize}
\item If $0<\alpha<2$, then $q_T<\infty$ precisely when
\begin{align}\label{multdimension}
 d<\frac{2\beta(2\alpha-1)}{\alpha}.
\end{align}
 \textcolor{blue}{In particular, the existence of a mild solution requires
\begin{equation}
    d < \frac{2\beta(2\alpha-1)}{\alpha}.
    \label{eq:mild-condition}
\end{equation}
Consequently:
    \item For $\alpha \leq \frac{1}{2}$, the inequality cannot be satisfied 
    for any $d \in \mathbb{N}$; hence the solution is  not mild.
    \item For $\alpha = 1$ and $\beta = 1$, the condition reduces to $d < 2$, 
    so the mild solution exists only in dimension $d = 1$.}
\item If $\alpha=2$, then $q_T<\infty$ precisely when $d<2\beta$.  Thus for
the classical Laplacian, $\beta=1$, the space--time white-noise wave equation
has a random-field mild solution only in dimension $d=1$.
\end{itemize}
\end{theorem}

\dproof
For $0<\alpha<2$, standard Mittag--Leffler estimates (see \cite{J}) and the change of variables
$\eta=s^{\alpha/(2\beta)}\xi$ in \eqref{multdalang} gives
\begin{align}
 \int_{\R^d}|\widehat K_\alpha(s,\xi)|^2d\xi
 =C_{\alpha,\beta,d,\lambda}
 s^{2\alpha-2-\alpha d/(2\beta)}.
\end{align}
The time integral near zero is finite  under
\eqref{multdimension}.  At $\alpha=2$,
\begin{align}
 \widehat K_2(s,\xi)
 =\frac{\sin(\sqrt\lambda |\xi|^\beta s)}
 {\sqrt\lambda |\xi|^\beta},
\end{align}
whose square is integrable in $\xi$ exactly when $d<2\beta$.
\fproof

\subsection{When $\alpha=2$}\label{alpha2}
Let $\beta=1$, $b=0$, and $c=\sqrt\lambda$.  Then
\begin{align}
 E_{2,1}(-c^2|\xi|^2t^2)&=\cos(c|\xi|t),\nonumber\\
 tE_{2,2}(-c^2|\xi|^2t^2)
 &=\frac{\sin(c|\xi|t)}{c|\xi|}.
\end{align}
Consequently equation \eqref{heat1} becomes
\begin{align}\label{waveeq}
 \partial_t^2Y(t,x)=c^2\Delta Y(t,x)
 +\sigma Y(t,x)\diamond W(t,x).
\end{align}
In dimension one,
\begin{align}
 K_2(t,x)=\frac{1}{2c}\mathbf{1}_{\{|x|<ct\}},
\end{align}
and the mild equation takes the explicit light-cone form
\begin{align}\label{wavecone}
 Y(t,x)&=\frac{u_0(x-ct)+u_0(x+ct)}{2}
 +\frac{1}{2c}\int_{x-ct}^{x+ct}u_1(z)\,dz\nonumber\\
 &\quad+\frac{\sigma}{2c}\int_0^t
 \int_{x-c(t-s)}^{x+c(t-s)}Y(s,z)\,B(ds,dz).
\end{align}
\textcolor{blue}{This formula exhibits finite propagation speed.  If $u_0=\delta_0$, the
deterministic part consists of two moving Dirac fronts; hence it is naturally
distribution valued.  A random-field model should use a regular initial pulse.}

 \section{Acknowledgments}
\begin{itemize}
\item
    \textcolor{blue}{
    We are grateful to Kenneth Karlsen for helpful suggestions and valuable comments.}
    \item 
 This research was funded by the Research Council of Norway under project 363480 (NASTRAN-Talent).
 \end{itemize}

\end{document}